\documentclass[11pt,leqno]{article}

\usepackage{amsmath,amsfonts,amscd,amssymb,theorem}

\long\def\comment#1\endcomment{}

\comment
\endcomment

\makeatletter
\begingroup
\gdef\th@dotted{\normalfont\itshape
  \def\@begintheorem##1##2{%
        \item[\hskip\labelsep \theorem@headerfont ##1\ ##2.]}%
\def\@opargbegintheorem##1##2##3{%
   \item[\hskip\labelsep \theorem@headerfont ##1\ ##2\ (##3).]}}
\endgroup
\makeatother

\theoremstyle{dotted}

\newtheorem{theorem}{Theorem}[section]
\newtheorem{lemma}[theorem]{Lemma}

\newtheorem{prop}[theorem]{Proposition}

\makeatletter
\begingroup
\gdef\th@upshape{\normalfont
  \def\@begintheorem##1##2{%
        \item[\hskip\labelsep \theorem@headerfont ##1\ ##2.]}%
\def\@opargbegintheorem##1##2##3{%
   \item[\hskip\labelsep \theorem@headerfont ##1\ ##2\ (##3).]}}
\endgroup
\makeatother

\theoremstyle{upshape}

\newtheorem{defn}[theorem]{Definition}
\newtheorem{remark}[theorem]{Remark}

\makeatletter
\renewcommand{\subsection}{\@startsection{subsection}{2}{0pt}{-3ex
plus -1ex minus -0.2ex}{-2mm plus -0pt minus
-2pt}{\normalfont\bfseries}} \makeatother

\makeatletter
\@addtoreset{equation}{section}
\makeatother

\newcommand{\cntrct}                % contraction with a vector field
{\hspace{2pt}\raisebox{1pt}{\text{$\lrcorner$}}\hspace{2pt}}

\newcommand{\proof}[1][Proof.]{\smallskip\noindent{\em #1}}
\def\endproof{\hfill\ensuremath{\square}\par\medskip}

\renewcommand{\labelenumi}{{\normalfont(\roman{enumi})}}
\renewcommand{\theenumi}{{\normalfont(\roman{enumi})}}

\def\eqref#1{\thetag{\ref{#1}}}

\let\latexref=\ref
\def\ref#1{{\normalfont{\latexref{#1}}}}

\newcommand{\wt}{\widetilde}

\newcommand{\hdot}{{\:\raisebox{3pt}{\text{\circle*{1.5}}}}}
\newcommand{\C}{{\mathbb C}}

\newcommand{\Q}{{\mathbb Q}}

\newcommand{\B}{{\cal B}}
\newcommand{\A}{{\cal A}}

\newcommand{\N}{{\cal N}}
\newcommand{\M}{{\cal M}}
\newcommand{\LL}{{\cal L}}
\newcommand{\ZZ}{{\cal Z}}

\newcommand{\X}{{\mathfrak X}}

\newcommand{\I}{{\cal I}}
\newcommand{\E}{{\cal E}}

\newcommand{\G}{{\sf G}}
\newcommand{\U}{U^0}

\newcommand{\calo}{{\cal O}}

\newcommand{\g}{{\mathfrak g}}

\newcommand{\gl}{{\mathfrak g}{\mathfrak l}_r}
\newcommand{\ssp}{{\mathfrak s}{\mathfrak p}}
\newcommand{\pglo}{{\mathfrak s}{\mathfrak l}}
\newcommand{\pgl}{\pglo_r}
\newcommand{\sltwo}{\pglo_2}

\newcommand{\z}{{\text{\bf z}}}
\newcommand{\w}{{\text{\bf w}}}

\renewcommand{\phi}{\varphi}

\renewcommand{\dim}{\operatorname{\sf dim}}
\newcommand{\codim}{\operatorname{\sf codim}}

\newcommand{\id}{\operatorname{\sf id}}
\newcommand{\rk}{\operatorname{\sf rk}}

\newcommand{\Hom}{\operatorname{Hom}}
\newcommand{\Ext}{\operatorname{Ext}}
\newcommand{\ext}{\operatorname{{\cal E}{\it xt}}}
\newcommand{\End}{\operatorname{End}}
\newcommand{\Aut}{\operatorname{Aut}}
\newcommand{\Stab}{\operatorname{Stab}}
\newcommand{\RHom}{\operatorname{RHom}}
\newcommand{\Rhom}{\operatorname{{\bf R}{\cal H}{\it om}}}

\newcommand{\hh}{\operatorname{{\mathcal H}{\mathcal H}}}
\newcommand{\hhom}{\operatorname{{\cal H}{\it om}}}

\newcommand{\Ker}{\operatorname{{\sf Ker}}}
\renewcommand{\Im}{\operatorname{{\sf Im}}}

\newcommand{\Bl}{\operatorname{Bl}}
\newcommand{\git}{/\!\!/}

\title{Local structure of hyperk\"ahler singularities in O'Grady's
examples}

\author{D. Kaledin\thanks{Partially supported by CRDF grant
RUM1-2694-MO05.} and M. Lehn}

\date{{\em To Victor Ginzburg, on the occasion of his 50-th birthday}}

\begin{document}

\maketitle

\tableofcontents

\section{Introduction}

An {\em irreducible hyperk\"ahler} manifold is by definition a
compact Riemannian manifold of dimension $4n$ with holonomy group
$Sp(n)$. Manifolds of this type are known to possess many additional
structures. In particular, an irreducible hyperk\"ahler $X$ is
automatically complex and K\"ahler. Moreover, the space $H^{2,0}(X)$
of global holomorphic $2$-forms is $1$-dimensional and spanned by a
single non-degenerate symplectic form $\Omega$. Complex manifolds
with this property are called {\em irreducible holomorphically
symplectic}.

Perhaps not unexpectedly for special holonomy manifolds, compact
irreducible hyperk\"ahler manifolds are very hard to construct -- so
much so that a large part of the structure theory of hyperk\"ahler
manifolds was developed by F. Bogomolov in the course of proving
that none exist in complex dimensions $2n > 2$. 

Fortunately for those who prefer positive statements, the last part
of the paper \cite{bogo} contained a mistake. This was noticed
several years later by Bogomolov himself and by A. Beauville. In the
meanwhile, the celebrated S.-T. Yau's proof of the Calabi conjecture
has arrived, and special holonomy metrics became much more
accessible. In particular, by a well-known argument Yau's Theorem
implied that a compact K\"ahler holomorphically symplectic manifold
admits a hyperk\"ahler metric. This immediately gives two examples
in complex dimension $2$: a torus $T$ and a surface $S$ of type
$K3$.  However, in contrast to the related theory of Calabi-Yau
manifolds, higher-dimensional examples are still not easy to come by
-- the purely algebro-geometric problem of finding compact
holomorphically symplectic manifolds is itself very hard.

A. Beauville \cite{beau} has overcome the difficulties and
constructed two series of examples of irreducible hyperk\"ahler
manifolds, each containing one manifold in each dimension $2n$, $n
\geq 2$. To understand these examples, recall that for any complex
surface $S$, the Hilbert scheme $S^{[n]}$ of $0$-dimensional
subschemes of length $n$ is a smooth complex variety of complex
dimension $2n$. Beauville's examples are:
\begin{enumerate}
\renewcommand{\labelenumi}{(\alph{enumi})}
\renewcommand{\theenumi}{(\alph{enumi})}
\item\label{hlb.ex} The Hilbert scheme $S^{[n]}$ of a $K3$-surface
$S$, and
\item The so-called {\em generalized Kummer variety} $K^{[n]}$ --
  namely, the pre-image of $0 \in T$ under the natural summation map
  $\Sigma:T^{[n+1]} \to T$ from the Hilbert scheme $T^{[n]}$ of a
  $2$-dimensional complex torus $T$ to $T$ itself.
\end{enumerate}
Beauville has proved that for every $n \geq 2$, both these manifolds
are K\"ahler and irreducible holomorphically symplectic, hence
irreducible hyperk\"ahler.

\bigskip

For a long time, the Beauville examples and their deformations were
the only known examples of irreducible hyperk\"ahler manifolds. The
situation was rather frustrating and also a bit strange, because at
least the example \ref{hlb.ex} admits an obvious generalization.
Instead of the Hilbert scheme $S^{[n]}$ for a $K3$-surface $S$, one
can consider the moduli space $\M_S(r,c_1,c_2)$ of torsion-free
sheaves on $S$ of rank $r$, with Chern classes $c_1$ and
$c_2$. These spaces are naturally symplectic outside of their
singular loci. Thus to produce a hyperk\"ahler manifold, all one has
to do is to choose the topological type $(r,c_1,c_2)$ in such a way
that the singular locus is empty, and apply Yau's
Theorem. Unfortunately, all the hyperk\"ahler manifolds obtained by
this procedure are not new -- it turns out that they are
deformationally equivalent to an appropriately chosen Hilbert
scheme.

\bigskip

The great breakthrough was achived by K. O'Grady about five years
ago. His new approach consisted of taking a {\em singular} moduli
space of sheaves on a $K3$ surface and resolving the singularities
in an appropriate way. The moduli space he considered was
$\M_S(2,0,4)$, the space of torsion-free sheaves $\E$ of rank $2$ on
a $K3$-surface $S$, with $c_1(\E) = 0$ and $c_2(\E)=4$. This is a
singular complex variety of dimension $10$, and in fact, this is the
simplest topological type for which the moduli space is singular. In
two papers \cite[part I, II]{ogr1}  (see also the published version 
\cite{ogr3})
O'Grady was able to construct a smooth symplectic desingularization
of the space $\M_S(2,0,4)$ (part I) and to prove that the resulting
hyperk\"ahler manifold is irreducible and not deformationally
equivalent to a manifold of Beauville's (part II). Later on, he also
produced a generalization of the Kummer variety construction, thus
obtaining a new irreducible hyperk\"ahler manifold of complex
dimension $6$ (\cite{ogr2}).

These examples are rather unique. In both cases, the singularity one
has to resolve is the same. O'Grady also considered moduli spaces
$\M(2,0,2k)$ of sheaves with $c_2=2k$ for $k > 2$; however, he was
not able to resolve their singularities in a symplectic fashion, and
conjectured that such a resolution does not exist.

\bigskip

The present paper arouse out of our attempts to understand the
O'Grady construction and to study the spaces $\M_S(2,0,2k)$ for
higher values of $k$. We focus on the local structure of the
singularities. We give a relatively explicit description of the
singularities in terms of the nilpotent coadjoint orbits of the
group $Sp(k)$. Then we prove -- using an idea of O'Grady himself --
that he was right: the moduli space $\M_S(2,0,2k)$ does {\em not}
admit a smooth projective symplectic resolution for $k \geq 3$. We
also give some speculation as to what might be still possible to do
at least in the case $k=3$.

Our results are not really satisfactory because they are negative:
we would much prefer to construct a symplectic resolution.
Nevertheless, we believe that the present paper will be useful --
firstly, because it settles definitely the case of moduli of
rank-$2$ bundles on a $K3$-surface, for better or for worse, and
secondly, because it might serve as an introduction to the rather
technical papers \cite{ogr1}, \cite{ogr3}. We hope that this will
stimulate further research in the area. To us, it looks promising
and important enough to merit wider attention.

\subsection*{Note added in proof.}
For technical reasons, this paper has spent a rather long time in
the form of an electronic preprint.  After it had been first posted
to the web, we have learned that a completely different and very
successful approach to the existence of resolutions had been
developped by Jaeyoo Choy and Young-Hoon Kiem (\cite{CK1},
\cite{CK3}). Their method is based on motivic intergration, and
allows to give a different proof of the non-existence of
resolutions. Also, a general non-existence results for the
singularities of the moduli spaces of sheaves on a K3 surface has
been proved later in \cite{KLS}. However, some of our results, such
as the formality statement for the Hilbert scheme and an explicit
description of the O'Grady singularity are still unavailable in the
literature.

\section{Main results}

We will now state our results. Let $S$ be a projective algebraic
surface over $\C$ of type $K3$. Fix a polarization $H$ on $S$. For
every positive integer $k \geq 2$, let $\M_S(k)=\M_S(2,0,2k)$ be the
Gieseker compactification of the moduli space of vector bundles $\E$
on $S$ with $\rk \E = 2$, $c_1(\E)=0$ and $c_2(\E)=2k$. This is a
projective algebraic variety of dimension $\dim\M_S(k)=6k-2$. In
\cite{ogr1}, O'Grady has proved the following.

\begin{lemma}[{{\cite[Lemma]{ogr1}}}]\label{ogr.lemma}
Singular points $m \in \M_S(k)$ are in one-to one correspondence
with isomorphism classes of sheaves of the form
\begin{equation}\label{ogr.dec}
\E_m \cong \I_\z \oplus \I_\w,
\end{equation}
where $\I_\z$ and $\I_\w$ are ideal sheaves of $0$-dimensional
subschemes $\z,\w \subset S$ of length $k$.\endproof
\end{lemma}

Interchanging $\z$ and $\w$ in \eqref{ogr.dec} of course gives an
isomorphic sheaf. There are no other isomorphisms. Thus the singular
locus $\M_S^{sing}(k)$ is naturally identified with the symmetric
square $S^2(S^{[k]})$ of the Hilbert scheme of $k$ points on $S$.
This symmetric square is itself naturally stratified: we have the
image of the diagonal in $S^{[k]} \times S^{[k]}$ and the open
complement to this image. In all, the variety $\M_S(k)$ has a
stratification with three strata,
\begin{equation}\label{strata}
\M_S^{bad}(k) \subset \M_S^{sing}(k) \subset \M_S(k),
\end{equation}
with $\M_S^{bad}(k) \cong S^{[k]}$ and $\M_S^{sing}(k) \cong
S^2(S^{[k]})$. The complements $\M^{sing}_S(k) \setminus
\M_S^{bad}(k)$ and $\M_S(k) \setminus \M_S^{sing}(k)$ are
non-singular. 

One can characterize the stratification \eqref{strata} by
considering automorphisms of the corresponding sheaves. If a point
$m \in \M_S(k)$ lies outside of the singular locus $\M^{sing}_S(k)$,
the automorphism group $\Aut\E_m$ of the corresponding sheaf $\E_m$
consists only of scalars, $\Aut\E_m\cong \C^*$. For a point $m \in
\M^{sing}_S(k) \setminus \M_S^{bad}(k)$, the group $\Aut\E_m$ of the
sheaf $\E_m=\I_\z\oplus\I_w$ is a torus $\C^*\times\C^*$. Finally,
for a point $m \in \M_S^{bad}(k)$ in the most singular stratum, we
have $\Aut\E_m \cong GL(2,\C)$.

Take such a point $m \in \M_S^{bad}(k) \subset \M_S(k)$. Our main
result describes the structure of the variety $\M_S(k)$ in a small
neighborhood of the point $m$. To state it, consider a symplectic
vector space $V = \C^{2k}$ of dimension $2k$, and let $Sp(V)=Sp(k)$
be the corresponding symplectic group, with its Lie algebra
$\ssp(V)$. Identify $\ssp(V)\cong\ssp(V)^*$ by means of the Killing
form. Denote by $\N_k \subset \ssp(V)$ the set of all symplectic
endomorphisms $B:V \to V$ such that
\begin{equation}\label{B}
B^2 = 0.
\end{equation}
The set $\N_k$ is a closed $Sp(V)$-invariant subvariety in
$\ssp(V)$. The only invariant of a symplectic endomorphism $B$
satisfying \eqref{B} is its rank. Thus we have a decomposition
$$
\N_k = \coprod \N^p_k
$$
into coadjoint orbits 
$$
\N^p_k = \left\{ B \in \ssp(V) \mid B^2=0 \text{ and } \rk B = p
\right\} \subset \ssp(V)\cong \ssp(V)^* 
$$
of the group $Sp(V)$. The unions
$$
\N_k^{\leq p} = \bigcup_{q \leq p} \N_k^q \subset \N_k
$$
form a stratification of the variety $\N_k$ by closed
$Sp(V)$-invariant strata. Since $B^2 = 0$, the image $B(V) \subset
V$ of an endomorphism $B \in \N^p_k$ must be an isotropic subspace in
$V$. Therefore $\N^p_k$ is empty unless $k \geq p$. When $p=0$, the
orbit $\N^0_k$ consists of a single point $0 \in \ssp(V)^*$.

\begin{prop}\label{local.prop}
Let $U_m \subset \M_S(k)$ be a small analytic neighborhood of $m \in
\M_S^{bad}(k) \subset \M_S(k)$.

Then there exists an open neighborhood $U \subset \N_k^{\leq 3}
\times V$ of $0$ in $\N^{\leq 3}_k \times V$, the product of the
stratum $\N^{\leq 3}_k \subset \ssp(V)^*$ with the symplectic vector
space $V$, and a map $\mu:U_m \to U$ such that
\begin{enumerate}
\item The map $\mu:U_m \to U$ is a two-to-one \'etale cover over $U
\cap \left(\N^3_k \times V\right)$.
\item Over $U \cap \left(\N^{\leq 2}_k \times V\right)$, the map
$\mu:U_m \to U$ is one-to-one
\item The preimage $\mu^{-1}\left(\N^1_k \times V\right) \subset
U_m$ is the intersection $U_m \cap \M_S^{sing}(k)$, and the preimage
$\mu^{-1}(0 \times V) \subset U_m$ is the intersection $U_m \cap
\M_S^{bad}(k)$.
\end{enumerate}
Moreover, the map $\mu$ sends the natural symplectic form on the
smooth part on $U_m$ to the natural symplectic from on $\N_k \times
V$ obtained from the given form on $V$ and the Kostant-Kirillov
symplectic form on the coadjoint orbit.
\end{prop}

The statement of the Proposition hides an essential difference
between the cases $k=2$ and $k \geq 3$. Namely, as noted above, in
the case $k=2$ the orbit $\N^3_k \subset \ssp(V)^*$ is
empty. Therefore $\N^{\leq 2}_k \times V = \N^{\leq 3}_k \times V$,
and we have a map $\mu:U_m \to \N^{\leq 2}_k \times V$ which is
generically one-to-one. The variety $\N^{\leq 2}_k$ has a natural
smooth symplectic resolution -- namely, the total space $T^*\G(V)$
of the Grassmanian $\G(V)$ of Lagrangian subspaces $L \subset V$ in
the symplectic vector space $V = \C^4$. This is essentially the
resolution constructed by O'Grady, although he did not use the
nilpotent orbit interpretation.

Recently this resolution also appeared in \cite{fu} as a part of a
beautiful general theorem: roughly speaking, for any simple
algebraic group $G$, {\em all} projective symplectic resolutions of
nilpotent coadjoint $G$-orbits are of the form $T^*(G/P)$ for some
parabolic subgroup $P \subset G$.

In the case $k \geq 3$, Proposition~\ref{local.prop} only gives a
map $\mu:U_m \to \N^{\leq 3} \times V$ which is generically
two-to-one. In this case we have the following.

\begin{theorem}\label{main}
For $k \geq 3$, the two-fold cover $U_m$ of the open neighborhood $U
\subset \N^{\leq 3}_k \times V$ considered in
Proposition~\ref{local.prop} does not admit a smooth projective
symplectic resolution compatible with the given symplectic
form. Consequently, for $k \geq 3$ the moduli space $\M_S(k)$ does
not admit a smooth symplectic desingularization $\wt{\M_S(k)}$.
\end{theorem}

We note that in the case $k=3$, the variety $\N^{\leq 3}_k$ itself
does admit a smooth symplectic desingularization -- again, it is the
cotangent space $T^*\G(V)$ of the Lagrangian Grassmannian
$\G(V)$. However, passing to a $2$-fold cover introduces
ramification, and the symplectic form is no longer
non-degenerate. As we will see in the proof of Theorem~\ref{main},
this cannot be cured.

We will prove Proposititon~\ref{local.prop} in
Section~\ref{local.sec}, after analyzing the deformation theory of
singular sheaves on $S$ in Section~\ref{def.sec}. In the course of
the analysis, we need a technical result on formality
(Propositon~\ref{frm}) whose proof we postpone so as not to
interrupt the expostion. The proof of Theorem~\ref{main} is
contained in the Section~\ref{main.sec}. The last
Section~\ref{frm.sec} is taken up with the proof of
Proposition~\ref{frm}.

\section{Local deformations of the most degenerate sheaf}\label{def.sec}

Fix a point $m \in \M_S^{bad}(k)$ in the most singular stratum in
the moduli space $\M_S(k)$, and let $\E_m$ be the corresponding
sheaf. Our first task is to analyze the deformation theory of the
sheaf $\E_m$.

By Lemma~\ref{ogr.lemma}, we have $\E_m \cong \I_\z \oplus \I_\z
\cong \I_\z \otimes \C^2$, where $\I_\z$ is the ideal sheaf of some
$0$-dimensional subscheme $\z \subset S$ of length $k$. Throughout
this section, we will work in wider generality and assume that $\E_m
= \I_\z \otimes \C^r$ for some fixed integer $r \geq 2$. It will 
affect neither the results nor the proofs.

\subsection{Generalities on deformations.}\label{gen.def}
We recall relevant general results from the deformation theory of
sheaves. Every sheaf $\E$ on a smooth projective complex manifold
admits a local deformation space $U_\E$ with a fixed origin point $o
\in U_\E$; the pair $\langle U_\E,o\rangle$ parametrizes small
deformations of the sheaf $\E$. To construct space $U_\E$, one
starts with the complex vector space $\Ext^1(\E,\E)$. By a standard
algebraic procedure, the DG Lie algebra structure on the complex
$\RHom^\hdot(\E,\E)$ defines a certain closed affine subvariety
$\wt{U}_\E \subset U_0$ in a small neighborhood $U_0 \subset
\Ext^1(\E,\E)$ of $0 \in \Ext^1(\E,\E)$. The variety $\wt{U}_\E$ is
called {\em a versal deformation space} of the sheaf $\E$. The
automorphisms group $G=\Aut(\E)$ of the sheaf $\E$ acts on the Lie
algebra $\Ext^\hdot(\E,\E)$ preserving the versal deformation space
$\wt{U}_\E \subset \Ext^1(\E,\E)$. The local deformation space
$U_\E$ is the affine algebraic variety obtained as the quotient
$\wt{U}_\E/G$. The origin point $o \in U_\E$ is the image of $0 \in
\Ext^1(\E,\E)$.

Given a point $u \in U_\E$ in the local deformation space, one can
also recover the automorphism group $\Aut\E_u$ of the corresponding
deformed sheaf $\E_u$ from the general formalism. To do this, one
notes that the $\Aut(E)$-action on the variety $\wt{U}_\E$ is not
free. Choose a lifting $\wt{u} \in \wt{U}_\E$ of the point $u \in
U_\E$. The stabilizer group $\Stab(\wt{u}) \subset \Aut(\E)$ does
not depend on the chosen lifting $\wt{u}$ and coincides with
$\Aut(\E_u)$ (while the embedding $\Aut(\E_u) \cong \Stab(\wt{u})
\subset \Aut(\E)$ does depend on the choice of $\wt{u}$ and is
determined by $u$ only up to a conjugation).

The algebraic procedure which gives a versal deformation space
$\wt{U}_\E$ is particularly simple in the case when the DG Lie
algebra $\RHom^\hdot(\E,\E)$ is {\em formal} -- in other words,
quasiisomorphic to its homology Lie algebra $\Ext^\hdot(\E,\E)$. In
this case, the Yoneda commutator bracket in $\Ext^\hdot(\E,\E)$
gives a quadratic map $Q:\Ext^1(\E,\E) \to \Ext^2(\E,\E)$. The zero
locus $Q^{-1}(0)$ of the map $Q$ is a versal deformation space
$\wt{U}_\E$ for the sheaf $\E$. In particular, $\wt{U}_\E \subset
\Ext^1(\E,\E)$ is a conic subvariety in the whole space
$\Ext^1(\E,\E)$, not only in a small neighborhood of $0$. (However,
the quotient $U_\E = \wt{U}_\E/G$ parametrizes deformations only in
a small neighorhood of the origin point.)

\subsection{The $\Ext$-algebra of the most degenerate sheaf.} 
We apply the general theory to a sheaf $\E=\E_m$ of the form $\E_m
\cong \I_\z \otimes \C^r$. Our first objective is to describe the DG
Lie algebra $\RHom^\hdot(\E_m,\E_m)$. First, consider
$\RHom^\hdot(\E_m,\E_m)$ as an associative DG algebra. We have
$$
\RHom^\hdot(\E_m,\E_m) \cong \End(\C^r) \otimes
\RHom^\hdot(\I_\z,\I_\z). 
$$
It turns out that in order to study the deformations of $\E_m$, all
we need to know is the $\Ext$-algebras $\Ext^\hdot(\E_m,\E_m)$ and
$\Ext^\hdot(\I_\z,\I_\z)$. This follows from the following general
fact.

\begin{prop}\label{frm}
Let $\I_\z$ be the ideal sheaf of a subscheme $\z \subset S$ of
length $0$ in a $K3$-surface $S$. Then the associative DG algebra
$\RHom^\hdot(\I_\z,\I_\z)$ is formal.
\end{prop}

The proof of this Proposition is slightly more technical than one
would like. We will postpone it till the last Section, so as not to
interrupt the exposition.

Granted Proposition~\ref{frm}, we wee that the DG algebra
$\RHom^\hdot(\E_m\E_m)$ is also formal and quasiisomorphic to the
product
\begin{equation}\label{exty}
\Ext^\hdot(\E_m,\E_m) \cong \End(\C^r) \otimes
\Ext^\hdot(\I_\z,\I_\z). 
\end{equation}
The $\Ext$-algebra $\Ext^\hdot(\I_\z,\I_\z)$ is well-known. We have
$\Ext^p(\I_\z,\I_\z) = 0$ unless $p=0,1,2$. For $p=0$ and $p=2$ we
have natural isomorphisms
$$
\Ext^0(\I_\z,\I_\z) \cong \Ext^2(\I_\z,\I_\z) \cong \C.
$$
The space $\Ext^1(\I_\z,\I_\z)$ is a complex vector space of
dimension $2k$. The Yoneda multiplication
$$
\Ext^1(\I_\z,\I_\z) \otimes \Ext^1(\I_\z,\I_\z) \to
\Ext^2(\I_\z,\I_\z) \cong \C
$$
is given by a non-degenerate symplectic form on the vector space
$V$. 

It will be convenient to introduce special notation for algebras of
this type. Let $V$ be a symplectic vector space. Define a
supercommutative algebra $\A^\hdot(V)$ by setting
$$
\A^p(V) \cong 
\begin{cases}
\C, &\qquad p = 0, \\
V, &\qquad p = 1, \\ 
\C, &\qquad p = 2,
\end{cases}
$$
with multiplication $\A^1(V) \otimes \A^1(V) = V \otimes V \to
\A^2(V) = \C$ obtained from the symplectic form on $V$. 

With this notation, we have $\Ext^\hdot(\I_\z,\I_\z) \cong
\A^\hdot(\Ext^1(\I_\z,\I_\z))$. Since this algebra is
supercommutative, the supercommutator bracket in the algebra
$\Ext^\hdot(\E_m,\E_m)$ is the product of the Lie bracket on the Lie
algebra $\End(\C^r) \cong \gl(\C)$ and the Yoneda multiplication in
$\Ext^\hdot(\I_\z,\I_\z)$.

It is convenient and customary to use the trace map $\gl(\C) \to \C$
to decompose the Lie algebra $\Ext^\hdot(\E_m,\E_m)$ into the sum of
the abelian Lie algebra $\Ext^\hdot(\I_\z,\I_\z)$ and the traceless
part $\Ext_0^\hdot(\E_m,\E_m)$. The summand $\Ext^\hdot(\I_\z,\I_\z)$
controls deformations $\wt{\I_\z}$ of the sheaf $\I_\z$ and
deformations of the sheaf $\E_m \cong \I_\z \otimes \C^r$ which are
of the form $\wt{\E_m} \cong \wt{\I_\z} \otimes \C^r$. The trace
decomposition is compatible with the Lie bracket and induces a
decompositon $U_m \cong U_0 \times \U_m$, where $U_0$ is the local
deformation space of the sheaf $\I_\z$, and $\U_m$ is the
deformation ``in transversal directions'' to $U_0 \subset
U_m$. Since the Hilbert scheme $S^{[k]}$ is smooth, the space $U_0$
is simply a neighborhood of $0$ in the vector space
$\Ext^1(\I_\z,\I_\z)$. All the information about singularities is
contained in the transversal deformation space $\U_m$, obtained from
the Lie algebra $\Ext_0^\hdot(\E_m,\E_m)$.

\subsection{Hamiltonian interpretation.}
By virtue of Proposition~\ref{frm}, in order to study the
deformation theory of the sheaf $\E_m$, it suffices to describe the
Yoneda square map in the traceless Lie algebra
$\Ext_0^\hdot(\E_m,\E_m)$. To do this, it is convenient to use the
language of Hamiltonian group actions.

Let $G$ be a reductive Lie group with Lie algebra $\g$, and let $V$
be a symplectic vector space. Identify $\g \cong \g^*$ by means of
the Killing form. Let $G$ act on the vector space $\g \otimes V$ by
the product of the adjoint representation and the trivial
one. Consider the symplectic form on $\g \otimes V$ obtained as the
product of the Killing form on $\g$ and the given symplectic form on
$V$. Then the action of $G$ on $\g \otimes V$ is Hamiltonian, and we
have a quadratic {\em moment map} $Q:\g \otimes V \to \g \cong
\g^*$.

On the other hand, define a graded Lie algebra $\LL^\hdot(\g,V)$ by
setting
\begin{equation}\label{L(g,V)}
\LL^\hdot(\g,V) = \g \otimes \A^\hdot(V).
\end{equation}

\begin{lemma}\label{yoneda.moment}
The moment map $Q:\g \otimes V \to \g$ for the $G$-action on $\g
\otimes V$ coincides with the Yoneda square map in the Lie algebra
$\LL^\hdot(\g,V)$.
\end{lemma}

\proof{} Clear. \endproof

Recall that whenever one has a Hamiltonian action of a reductive Lie
group $G$ on an affine symplectic manifold $X$, one defines the {\em
Hamiltonian reduction} $X\git{}G$ as the quotient of the preimage
$Q^{-1}(0) \subset X$ of $0$ under the moment map $Q:X \to \g^*$ by
the action of $G$. Denote by $P(G,V) = \g \otimes V \git{}G$ the
Hamiltonian reduction of the space $\g \otimes V$ by the $G$-action.

By \eqref{exty}, the traceless part $\Ext_0^\hdot(\E_m,\E_m)$ of the
Lie algebra $\Ext^\hdot(\E_m,\E_m)$ is isomorphic to
$\LL^\hdot(\pgl,V)$, $V = \Ext^1(\I_\z,\I_\z)$. The automorphism
group $\Aut(\E_m)\cong GL(n)$ acts on $\LL^\hdot(\pgl,V)$ through
its quotient $G=PGL(n)$. Thus we are in the situation described by
Lemma~\ref{yoneda.moment}. We conclude that the Yoneda bracket map
$Q$ for the Lie algebra $\LL^\hdot(\pgl,V)$ coincides with the
moment map for the $G$-action on $\pgl \otimes V$. Since the DG Lie
algebra $\RHom^\hdot(\E_m,\E_m)$ is formal, this implies that the
transversal local deformation space $\U_m$ of the sheaf $\E_m$ is
canonically an open neighborhood $\U_n \subset P(G,V)$ of $0$ in the
Hamiltonian reduction $P(G,V)$. To sum up, we have proved the
following.

\begin{prop}\label{=redukcii}
The local deformation space $U_m$ of the sheaf $\E_m \cong \I_\z
\otimes \C^r$ is isomorphic to an open neighborhood of $0$ is the
product $P(PGL(r),V) \times V$ of the Hamiltonian reduction
$P(PGL(r),V)$ with the symplectic vector space $V$.\endproof
\end{prop}

We note that the product $P(PGL(r),V) \times V$ is itself a
Hamiltonian reduction: we have $P(PGL(r),V) \times V \cong
P(GL(r),V)$. It is also possible to prove that the embedding $U_m
\subset P(PGL(r),V) \times V \cong P(GL(r),V)$ is symplectic outside
of singularities. To do this, one has to unwind the definition of
the symplectic form on the deformation space $U_m$ and see that it
is the same form as the one obtained from the reduction. We leave
the details to the reader.

\begin{remark}
When $\dim V \geq 4$, one can also find the dimension of the space
$P(PGL(r),V)$ -- hence also of the deformation space $\U_m$. It is equal
to
$$
\dim \U_m = \dim P(PGL(r),V) = \dim \pgl \cdot (\dim V - 2).
$$
Indeed, decompose symplectically $V = \C^2 \oplus V'$. Since the Lie
algebra $\pgl$ is semisimple, we have $[\pgl,\pgl]=\pgl$. Thus for a
generic vector $v' \in V' \otimes \pgl$, one can find $a,b \in \pgl$
with $[a,b] = Q(v') \in \pgl$. Then $\langle a,b \rangle \oplus v'
\in V \otimes \pgl$ satisfies $Q(v) = 0$.  On the other hand, since
$v'$ is generic, the stabilizer $\Stab(v) \subset \Stab(v') \subset
PGL(r)$ of the vector $v$ is the trivial subgroup. Therefore the
differential of the moment map $Q:\pgl \otimes V \to \pgl$ is
surjective in the generic point of the zero locus $Q^{-1}(0)$, the
$PGL(r)$-action is generically free, and one has
$$
\begin{aligned}
\dim P(PGL(r),V) &= \dim \pgl \otimes V - \dim PGL(r) - \dim \pgl \\
&= \dim \pgl \otimes V - 2\dim \pgl.
\end{aligned}
$$
This coincides with the expected dimension of the moduli space
computed by Riemann-Roch. It is interesting to note that when $\dim
V = 2$, the dimension of the space $P(PGL(r),V)$ is actually
strictly bigger than the right-hand side of this equation (which in
this case is equal to $0$). Geometrically, this means that the sheaf
$\E_m$ for $k=1$ cannot be deformed to a non-singular vector bundle.
\end{remark}

\section{Hamiltonian action of the symplectic group}\label{local.sec}

By virtue of Proposition~\ref{=redukcii}, in order to study the
local deformation space of a sheaf $\E_m \cong \I_\z \otimes \C^r$
and to prove Proposition~\ref{local.prop}, it is sufficient to study
the Hamiltonian reduction space $P(G,V)$ with
$V=\Ext^1(\I_\z,\I_\z)$ and $G=PGL(r)$. This is the goal of this
section.

The symplectic group $Sp(V)$ acts naturally on the vector space $V$
preserving the symplectic form $\Omega$. This action is Hamiltonian:
the moment map $\mu_0:V \to \ssp(V)^*$ is quadratic and comes from
the standard identification $S^2V \cong \ssp(V) \cong
\ssp(V)^*$. Tensoring with the Lie algebra $\pgl$ gives a
Hamiltonian action of the group $Sp(V)$ on the symplectic vector
space $\pgl \otimes V$. The corresponding moment map $\mu:\pgl
\otimes V \to \ssp(V)^*$ is also quadratic and coincides with the
product of the moment map $\mu_0:V \to \ssp(V)^*$ with the Killing
form on $\pgl$. Moreover, the $Sp(V)$-action on $\pgl \otimes V$
commutes with the $PGL(r)$-action -- $Sp(V)$ acts on the second factor,
while $PGL(r)$ acts on the first one.

Thus we have a symplectic vector space $\pgl \otimes V$ with two
commuting Hamiltonian group actions -- firstly, a $PGL(r)$-action
with the moment map $Q:\pgl \otimes V \to \pgl$, secondly, a
$Sp(V)$-action with the moment map $\mu:\pgl \otimes V \to
\ssp(V)^*$. Our strategy is to study the reduction $P(G,V) = \pgl
\otimes V \git{}PGL(r)$ by means of the second moment map $\mu:\pgl
\otimes V \to \ssp(V)$.

From now on, we restrict our attention to the case of rank
$r=2$. Denote $G=PGL(2)$, $\g=\sltwo$. We will use the
identification $PGL(2) \cong SO(3)$ and think of $G$ as the group of
automorphisms of a $3$-dimensional complex vector space $W$ which
preserve a non-degenerate symmetric form $h \in S^2W^*$ on $W$ and
an orientation -- or, equivalently, an isomorphism $\chi:\Lambda^2W
\cong W$. One can think of $W$ as the Lie algebra $\g$, with its
Killing form and the commutator map $[-,-]:\Lambda^2\g\cong\g$.

Identify the symplectic vector space $\g \otimes V \cong W \otimes
V$ with the space
$$
\Hom(W^*,V)
$$
of linear maps from the dual space $W^*$ to $V$. In this
interpretation, it is very easy to describe the moment maps $Q$ and
$\mu$.

\begin{lemma}\label{moments}
Let $a \in \Hom(W^*,V) \cong \g \otimes V$ be an arbitrary map. Then
$$
Q(a) = \chi^*(a^*\Omega) \in W \cong \g,
$$
while
\begin{equation}\label{mu}
\mu(a) = a(h) \in S^2V \cong \ssp(V).
\end{equation}
\end{lemma}

\proof{} Clear. \endproof

We see that a map $a \in \Hom(W^*,V)$ lies in the zero set of the
moment map $Q$ if and only if the image $a(W^*) \subset V$ is an
isotropic subspace. By definition, the reduction space $P(G,V) =
\pgl \otimes V \git{}PGL(r)$ is the geometric invariant theory quotient
of the zero set $Q^{-1}(0)$ by the group $G$. In particular, this
means that points in $P(G,V)$ are in one-to-one correspondence with
the closed $G$-orbits in $Q^{-1}(0)$. In order to analyze these
orbits, we need a preliminary linear algebra result.

\begin{lemma}\label{closed}
  Let $U$, $U'$ be two finite-dimensional complex vector spaces that
  are equip\-ped with symmetric pairings. Assume that $\dim U \leq \dim
  U'$ and that the pairing on $U'$ is non-degenerate. Denote by $U_0
  \subset U$ the radical of the pairing on $U$.  Consider the space
  $O(U,U')$ of maps $f:U \to U'$ compatible with the pairings, and
  let the orthogonal group $O(U')$ act on $O(U,U')$ by
  multiplication on the right.
  
  Then $O(U,U')$ contains exactly one closed $O(U')$-orbit, and it
  coincides with the subset $O(U/U_0,U') \subset O(U,U')$ of maps $f
  \in O(U,U')$ which vanish on $U_0 \subset U$.
\end{lemma}

\proof{} This is a reasonably standard fact about parabolic
subgroups in the orthogonal group; however, we include a direct
linear-algebraic proof for the convenience of the reader. First of
all, it is obvious that $O(U)$ acts transitively on the set
$O(U/U_0,U')$, so that it indeed forms a single $O(U)$-orbit. A map
$f \in O(U,U')$ lies in $O(U/U_0,U')$ if and only if $\dim \Im f =
\dim(U/U_0)$. Since for any $f \in O(U,U')$ we must have $\dim \Im f
\geq \dim (U/U_0)$, the subset $O(U/U_0,U') \subset O(U,U')$ is
indeed closed. 

It remains to show that all the other $O(U)$-orbits in $O(U,U')$ are
{\em not} closed. Fix an element $f \in O(U,U')$ which does not lie
i $O(U/U_0,U') \subset O(U,U')$, and choose a complement $U_1
\subset U$ to the subspace $U_0 \subset U$. Then the pairing is
non-degenerate on $U_1$, so that the map $f:U_1 \to U'$ is
injective.  Let $U'_0 \subset U'$ be the orthogonal complement to
$f(U_1) \subset U'$. Then $f$ induces a non-trivial map $f_0:U_0 \to
U_0'$, we have a natural group embedding $O(U_0') \subset O(U')$,
and it suffices to check that $O(U'_0) \cdot f_0$ is not closed in
$O(U_0,U_0')$. In other words, we can assume that the pairing on $U$
is trivial (and $f \neq 0$).  Then the claim becomes obvious:
indeed, for any $f \in O(U,U')$ the group $O(U)$ contains a
one-parameter subgroup $e:\C^* \to O(U)$ such that $e(\lambda) f =
\lambda f$ for any $\lambda \in \C^*$.
\endproof

\begin{lemma}\label{fibers}
  Let $a \in \Hom(W^*,V)$ lie in the zero level of the moment map
  $Q$, and assume that the orbit $O(W) \cdot a$ is closed in
  $Q^{-1}(0)$. Then
  $$
  \mu(a) \in \N^p_k \subset \ssp(V),
  $$
  where $p=\rk W'$ is the dimension of the image $W' = \Im a
  \subset V$. Moreover, $\Im \mu(a) = W'$ and the associated
  symmetric form $h(\mu(a)) \in S^2W'$ coincides with
  $a(h)$.
\end{lemma}

\proof{} For every $a \in \Hom(W^*,V)$ with image $W' = \Im a
\subset V$, \eqref{mu} shows that the moment value
$\mu(a)\in\ssp(V)$ considered as a symplectic automorphism $\mu(a):V
\to V$ satisfies $\Im\mu(a) \subset W'$. Equip $W'$ with the
symmetric form $h(\mu(a)) \in S^2W'$.  By Lemma~\ref{closed} (or
rather, by the dual statement), the orbit $O(W) \cdot a$ is closed in
$O(W,W')$ if and only if we have an equality, $\Im\mu(a) \subset
W'$.  This implies $\Ker\mu(a) = {W'}^\perp$, the symplectic
orthogonal to the subspace $W' \subset V$. If we have $Q(a)=0$, then
$W'$ is isotropic; therefore $W' \subset {W'}^\perp$ and
$\mu(a)^2=0$, so that $\mu(a) \in \N_k$.  Conversely, every
symplectic endomorphism $B \in \N_k$ with $W' = \Im B \subset V$
corresponds to a non-degenerate symmetric $h(B) \in S^2W' \subset
S^2V$ on the dual space ${W'}^*$ under the identification $S^2V
\cong \ssp(V)$.
\endproof

Thus the moment map $\mu$ restricts to a map $\mu:Q^{-1}(0) \to \N_k
\subset \ssp(V)$. Since it is $G$-invariant, it descends to the
reduction $P(G,V)$. 

\begin{lemma}
Let $B \in \N_k \subset \ssp(V)$ be a symplectic endomorphims $B:V
\to V$ satisfying $B^2=0$, let $W' = \Im B \subset V$ be its image,
and let $h(B) \in S^2W'$ be the non-degenerate symmetric form on
${W'}^*$ associated to $B$.

Then the fiber $\mu^{-1}(B) \subset P(G,V)$ is isomorphic to the set
$O(W',W)$ of embeddings ${W'}^* \hookrightarrow W$ with $a^*(h) =
h(B)$, modulo the natural action of $G=SO(W)$.
\end{lemma}

\proof{} By \eqref{mu}, every $a \in Q^{-1}(0) \subset \Hom(W^*,V)$
with $\mu(a) = B$ must send $h \in S^2W^*$ to $h(B)$ on $W' = \Im a
= \Im B$. In other words, the transposed map $a^*:{W'}^*
\hookrightarrow W$ satisfies $a^*(h)=h(B)$. Conversely, since
$B^2=0$, $W'$ is isotropic. Therefore by Lemma~\ref{moments} every
$a \in \Hom(W^*,V)$ with $\Im a = W'$ and $a^*(h) = h(B)$ satisfies
$Q(a)=0$ and $\mu(a)=B$. \endproof

We can now prove Proposition~\ref{local.prop}

\proof[Proof of Proposition~\ref{local.prop}.]  Embed $U_m$ into
$P(G,V) \times V$ by Proposition~\ref{=redukcii}. Take $B \in \N_k
\subset \ssp(V)$ with $W' = \Im B$.  By Lemma~\ref{fibers}, the
fiber of the map $\mu:P(G,V) \to \ssp(V)$ over the point $B$
coincides with the set $O(W',W)$ of Euclidean embeddings ${W'}^*
\hookrightarrow W$ modulo the natural $SO(W)$-action.

If $B \in \N^3_k$, so that $\dim W'=3$, the quotient set
$O(W',W)/SO(W)$ consists of two points (corresponding to the choice
of the orientation induced on $W'$). This proves \thetag{i}. If $B
\in \N^{\leq 2}_k$, so that $\dim W' \leq 2$, the action of $SO(W)$
on $O(W',W)$ is transitive, and this proves \thetag{ii}. However,
when $B \in \N^{\leq 1}_k$, any corresponding point $a^* \in
O(W',W)$ -- hence also any point $a \in Q^{-1}(0) \subset
\Hom(W^*,V)$ with $\mu(a) = B$ -- acquires a non-trivial stabilizer
in $G = PGL(2)$. If $\rk B = 1$, this stabilizer is $\C^* = \C^*
\times \C^*/\C^* \subset PGL(2) = GL(2)/\C^*$, while for $B = 0$ the
stabilizer is the whole $PGL(2)$. By Subsection~\ref{gen.def}, this
stabilizer $\Stab(a)$ is the image in $PGL(2)$ of the automorphism
group $\Aut(\E_{a,v}) \subset GL(2)$ of the sheaf $\E_{a,v}$
parametrized by the point $a \times v \in U_m \subset P(G,V) \times
V$ (no matter what is the vector $v \in V$). Hence
$\mu^{-1}(\N^1_k)$ parametrizes sheaves $\E$ with $\Aut\E = \C^*
\times \C^*$, and $\mu^{-1}(0)$ parametrizes sheaves with $\Aut\E =
GL(2)$. This proves \thetag{iii}.  \endproof

As a corollary, we see that the $Sp(V)$-action on $P(G,V)$ is
transitive. The two-to-one ramification in the moment map
$\mu:P(G,V) \to \N^{\leq 3}_k$ comes from the action of the element
$\iota = -\id:V \to V$, $\iota \in Sp(V)$. This element acts
trivially on the orbit $\N^3_k$ and non-trivially on the variety
$P(G,V)$.

\begin{remark}
As the reader can see, our proofs rely heavily on the identification
$PGL(2) \cong SO(3)$ and the isomorphism $\Lambda^2\sltwo \cong
\sltwo$. Thus they do not generalize to the case of higher rank
$r$. However, the main result, -- namely, the fact that the
$Sp(V)$-action on $P(G,V)$ is transitive, and the moment map
$\mu:P(G,V) \to \ssp(V)$ is generically finite onto its image, --
may generalize to the higher rank case, at least when the dimension
of the vector space $V$ is comparable to the dimension of the Lie
algebra $\pgl$. We conjecture that it is enough to require
$$
\dim V \geq 2 \dim \pgl = 2(r^2-1).
$$
This would give a cubic lower bound on the second Chern class $c_2 =
\frac{1}{2} r \dim V$.  

To be more precise, we note that the first part of our approach,
namely, Lemma~\ref{moments}, holds in full generality modulo an
appropriate change of notation. However, for $r \geq 3$ the
commutator map $\chi:\Lambda^2\pgl \to \pgl$ is no longer an
isomorphism. Therefore the condition $Q(a)=0$ no longer implies that
$\Im a \subset V$ is an isotropic subspace, but only imposes some
weaker condition on the subspace $\Im a \subset V$. Consequently,
the orbits one encounters in the image of the moment map $\mu$,
while still nilpotent, do not contain isomorphisms with $B^2=0$ but
something more complicated. Whether these orbits or their covers
admit symplectic resolutions should be the topic of further
research.
\end{remark}

\begin{remark} Proposition~\ref{=redukcii} is purely local. However,
we note that the O'Grady singularity $\N=\N^{\leq 2}_2 \subset
\ssp(\C^4)$ is $\Q$-factorial. Using this fact, it is easy to check
that the argument in \cite[Proposition 5.6]{K1} applies to $\N$ and
proves that the resolution $T^*\G(\C^4)$ is isomorphic
$$
T^*\G(\C^4) \cong \Bl(\N,\E_l)
$$
to the blow-up of the variety $\N$ in a certain sheaf of ideals
$\E_l \subset \calo_\N$. The sheaf $\E_l$ is obtained by taking an
appropriate positive integer $l$ and extending, by push-forward, the
$l$-th power of the ideal sheaf of the subvariety $\N^{\leq 1}_2
\subset \N$ from $\N \setminus \{0\}$ to $\N$. This shows how to
construct the global symplectic resolution of the moduli space
$\M_S(4)$: it suffices to take the $l$-th power of the ideal sheaf
of the subvariety
$$
\left( \M^{sing}_S(2) \setminus \M^{bad}_S(2) \right) 
\subset \left( \M_S(2) \setminus \M^{bad}_S(2) \right),
$$
extend it by push-forward from $\M_S(2) \setminus \M^{bad}_S(2)$ to
the whole $\M_S(2)$, and take the blow-up of the resulting ideal
sheaf. In this way one can avoid a reference to the Mori Cone
Theorem in \cite{ogr1}.
\end{remark}

\section{Existence of resolutions}\label{main.sec} 

Consider the Hamiltonian reduction space 
$$
U = P(GL(2),V) = P(PGL(2),V) \times V 
$$
which, by Proposition~\ref{local.prop}, gives a local model for the
singularity of the moduli space $\M(2,0,k)$. We will now prove
Theorem~\ref{main} claiming that for $k \geq 3$, the space $U$ does
not admit a smooth projective resolution compatible with the given
symplectic form.

We will need the following general facts on the geometry of a
projective symplectic resolution.

\begin{lemma}\label{vf.lifts}
Let $X$ be a smooth symplectic algebraic variety equipped with a
projective birational map $X \to Y$ to an irreducible normal affine
algebraic variety $Y$.
\begin{enumerate}
\item The map $X \to Y$ is semismall -- in other words, we have
$\dim X \times_Y X = \dim X$.
\item Every vector field $\xi$ on $Y$ lifts to a (unique) vector
field on $X$.
\end{enumerate}
\end{lemma}

\proof{} \thetag{i} is \cite[Proposition 1.2]{K}, and \thetag{ii} is
\cite[Lemma 5.3]{GK}. \endproof

This immediately implies the following. Assume given an algebraic
variety $Y$, and a non-degenerate symplectic form $\Omega$ on the
smooth locus $Y^{sm} \subset Y$; then we will say that a smooth
resolution $\pi:X \to Y$ is {\em compatible with the form $\Omega$}
if $\pi^*\Omega$ extends to a non-degenerate sympletic form on the
whole smooth algebraic variety $X$.

\begin{lemma}\label{factors}
Let $Y = Y' \times Z$ be the product of a smooth affine symplectic
variety $Z$ and an affine variety $Y'$ equipped with a symplectic
form on the smooth locus $Y'_{sm} \subset Y'$. Assume that $Y$
admits a projective smooth resolution $f:X \to Y$ compatible with
the product symplectic form on $Y'_{sm} \times Z$.

Then for an arbitrary point $z \in Z$, the preimage 
$$
X' = f^{-1}(Y' \times \{z\}) \subset X
$$
is a smooth projective resolution of the variety $Y'$ compatible
with the given symplectic form on $Y'_{sm} \subset Y'$.
\end{lemma}

\proof{} The only thing to prove is the smoothness of $X'$. In other
words, we have to prove that the composition $X \to Y \to Z$ of the
map $f:X \to Y$ with the natural projection $Y \to Z$ is a smooth
map over every point $z \in Z$. To prove it, choose a set of vector
fields $\xi_1,\ldots,\xi_l$ on $Y$ such that the vectors $\xi_i(z)$
form a basis of the tangent space $T_zZ$. By Lemma~\ref{vf.lifts},
all the vector fields $\xi_i$ lift to vector fields on
$X$. Therefore the differential of the projection $X \to Z$ is
surjective for every point $x \in X$ lying over $z \in Z$.
\endproof

Lemma~\ref{factors} immediately implies that to prove that $U =
P(PGL(2),V) \times V$ does not admit a symplectic resolution, it
suffices to prove that the first factor $P(PGL(2),V)$ does not admit
such a resolution. 

However, it turns out that we can do more -- namely, apply
Lemma~\ref{factors} to the space $P(PGL(2),V)$ itself. To do this,
denote $P=P(PGL(2),V)$ and let $P_1 = \mu^{-1}(\N^1_k) \subset P$ be
the singular stratum -- that is, the preimage of the stratum $\N^1_k
\subset \N^{\leq 3}_k$ under the $Sp(V)$-moment map $\mu:P \to
\N^{\leq 3}_k \subset \ssp(V)$. We want to show that locally near a
point $p \in P_1$, the space $P=P(PGL(2),V)$ admits a product
decomposition of the type required to apply Lemma~\ref{factors}.

Recall that in the language of Section~\ref{local.sec}, points in
the space $P$ correspond to maps $a:W^* \to V$ from the
Euclidean $3$-dimensional vector space $W^*$ to the symplectic
space $V$ with isotropic image and Euclidean kernel, 
considered modulo the action of the group $SO(W) \cong
PGL(2)$. Then a point $p \in P$ lies in the stratum $P_1$ if and
only if the corresponding subspace $a(W^*) \subset V$ is of
dimension exactly $1$. 

Fix once and for all a vector $w \in W^*$. We have the orthogonal
decomposition $W^* = \C \cdot w \oplus W^\perp$. For every point
$p \in P_1$, we can apply a suitable element $g \in SO(W)$ and
arrange so that the kernel $\Ker a$ of the corresponding map $a:W^*
\to V$ coincides with the $2$-dimensional subspace $W^\perp
\subset W^*$. Having done this, we are left with only one invariant
of the point $p \in P_1$ --namely, the vector $v_p = a(w) \in V$.
The stabilizer of the subspace $W^\perp$ in $SO(W)$ is the group
$O(W^\perp) \subset SO(W)$. It acts on the line $\C \cdot w
\subset W^*$ by multiplication by $\pm 1$, so that the vector $v_p
\subset V$ is well-defined up to a sign. Thus we have an unramified
two-fold cover $\wt{P}_1 \to P_1$ and an isomorphism $\wt{P}_1 \cong
V \setminus \{0\}$, $p \mapsto v_p$.

Let now $V'$ be a symplectic vector space of dimension $\dim V - 2 =
2(k-1)$, and consider the vector space $\Hom(W^\perp,V')$. Since
$W^\perp$ is Euclidean, $\Hom(W^\perp,V')$ is naturally a
symplectic vector space equipped with a Hamiltonian action of the
group $O(W^\perp)$. In particular, we have a Hamiltonian action of
the subgroup $\C^* \cong SO(W^\perp) \subset O(W^\perp)$. Denote
by $P' = \Hom(W^\perp,V')\git{}\C^*$ the variety obtained by
Hamiltonian reduction. As in Section~\ref{local.sec}, points $p' \in
P'$ correspond to maps $a':W^\perp \to V'$ with isotropic image.

\begin{lemma}\label{dec}
There exists a dense open subset $U \subset P$, a dense open subset
$U_1 \subset \wt{P}_1$ and an \'etale map
$$
\rho:U_1 \times P' \to U
$$
compatible with the natural symplectic forms on the smooth parts of
both sides.
\end{lemma}

\proof{} Fix a subspace $V_1 \subset V$ of codimension $1$. Its
symplectic orthogonal $V_1^\perp \subset V$ is a line in $V_1$, and
the quotient $V_1/V_1^\perp$ is a symplectic vector space of
dimension $\dim V - 2$. Fix an identification $V_1/V_1^\perp \cong
V'$. Let $U \subset P$ be the set of points $p \in P$ such that the
corresponding map $a:W^* \to V$ is transversal to $V_1$ -- in other
words, $V = a(W^*) + V_1$. This is true generically, so that $U
\subset P$ is a dense open subset. Set $U_1 = V \setminus V_1
\subset V \setminus \{0\} \cong \wt{P}_1$.

To define the map $\rho:U_1 \times P' \to U$, note that since for
every point $v_p \times p' \in U_1 \times P'$ the intersection $V_1
\cap v_p^\perp$ is transversal, so the projection induces a
canonical symplectic isomorphism $V_1 \cap v_p^\perp \cong V'$. In
other words, we have a canonical symplectic embedding $\chi_p:V' \to
v_p^\perp \subset V$ associated to the vector $v_p \in V \setminus
V_1$. If $a':W^\perp \to V'$ is the map associated to a point $p'
\in P'$, then the composition $\chi_p \circ a':W^\perp \to V_1$ is
isotropic embedding into $v_p^\perp \subset V$. We define $\rho(v_p
\times p')$ as the point $p \in P$ corresponding to the direct sum
map $(\chi_p \circ a') \oplus v_p:W^\perp \oplus \C\cdot w \to V$;
since $(\chi_p \circ a')(W^\perp)$ is isotropic and orthogonal to
$v_p$, the point $p \in P$ is well-defined.

We claim that the map $\rho$ is \'etale, moreover, it is an
unramified two-fold cover. Indeed, take a point $p \in P$, and let
$a:W^* \to V$ be the associated map. The point $p$ lies in the image
of the map $\rho$ if and only if $a(W^\perp)$ lies in the subspace
$V_1 \subset V$. But if $p \in U$, the preimage $a^{-1}(V_1) \subset
W^*$ is by definition a $2$-dimensional subspace. Therefore applying
a suitable $g \in SO(W)$, we can arrange so that $a^{-1}(V_1) =
W^\perp \subset W$, and modulo $SO(W^\perp)$, there are exactly two
choices of such a $g$.
\endproof

By Lemma~\ref{factors}, this Lemma implies that if the space
$P=P(PGL(2),V)$ were to admit a symplectic resolution, then so would
the space $P'$. Indeed, then the subset $U \subset P$ would also
have a symplectic resolution $X$, and the \'etale cover $\wt{X} = X
\times_U (U_1 \times P')$ would be a symplectic resolution for $U_1
\times P'$.

By itself, this is not enough to derive a contradiction and prove
Theorem~\ref{main} -- we need to throw in one more
ingredient. Namely, consider the element $\iota \in Sp(V)$ given by
$-\id:V \to V$. This element acts on $P$ interchanging the leaves of
the two-fold covering $\mu:P \to \N^{\leq 3}_k$. In particular, it
acts trivially on the points lying over $\N^{\leq 2}_k \subset
\N^{\leq 3}_k$, such as points $p \in P_1$ in the stratum $P_1 =
\mu^{-1}(\N^1_k) \subset P$.  Moreover, it is easy to see that
$\iota$ preserves $U \subset P$ and is compatible with the two-fold
\'etale cover $\rho:U_1 \times P \to U$ -- namely, there exists an
involution $\iota:P\ \to P'$ such that $\rho \circ (\id \times
\iota) = \iota \circ \rho$. To prove this, notice that $\iota$ acts
on the space
$$
P = (W \otimes V)\git{}PGL(2)
$$ 
by the endomorphism $-\id$ on $V$ or, which is equivalent, by the
endomorphism $-\id$ on $W$. Moreover, since we take a quotient by
the group $PGL(2)=SO(3)$, we can replace $-\id$ with an arbitrary
orthogonal map $\tau:W \to W$ with determinant $-1$. To get a model
compatible with the \'etale map $\rho$, we need to choose $\tau$ so
that it preserves $W^\perp \subset W^*$ and the orthogonal vector $w
\in W$. It suffices to take an arbitrary orthogonal map $W^\perp \to
W^\perp$ with determinant $-1$. The transposition $\tau$
interchanging $w_1$ and $w_2$ is a good choice. The involution
$\iota$ on $P' = \Hom(W^\perp,V')\git{}SO(W^\perp)$ is then induced
by composition with $\tau$ on $\Hom(W^\perp,V')$.

\begin{lemma}\label{no.res}
An arbitrarily small $\iota$-invariant neighborhood $U' \subset P'$
of $0 \subset P'$ does not admit a $\iota$-equivariant smooth
projective resolution $X' \to U'$ compatible with the given
symplectic form on the smooth locus $U'_{sm} \subset U'$.
\end{lemma}

\proof{} The space $W^\perp$ is a $2$-dimensional Euclidean vector
space over $\C$. The group $\C^*=SO(2)$ acts on $W^\perp$ with two
eigenvalues of opposite signs. Fix corresponding eigenvectors
$w_1,w_2 \in W$. The vectors $w_1,w_2$ have length $0$ and form a
basis of the vector space $W^\perp$. Therefore we have a
decomposition
$$
\Hom(W^\perp,V) \cong W^\perp \otimes V \cong w_1V \oplus w_2V,
$$
and both $w_1V$ and $w_2V$ are Lagrangian subspaces. Identifying $V
\cong V^*$ by means of the symplectic form, we obtain a
decomposition $W^\perp \otimes V \cong V \oplus V^*$. The symplectic
form in this interpretation is the standard form induced by the
duality between $V$ and $V^*$. It no longer depends on the
symplectic form on $V$. 

The Hamiltonian reduction $P' \cong (V \oplus V^*)/\C^*$ is
well-known -- it is a quadratic cone obtained by contracting the zero
section in the total space $T^*P(V)$ of the cotangent bundle to the
projectivization $P(V)$ of the complex vector space $V$. In
particular, the complement 
$$
P' \setminus \{0\} \cong T^*P(V) \setminus P(V)
$$ 
is smooth. Since $\codim P(V) \subset T^*P(V)$ is equal to $\dim
P(V) = 2k-1 \geq 2$, the Picard group of the complement $P'
\setminus \{0\}$ coincides with the Picard group of the projective
space $P(V)$. It is freely generated by the class of a single line
bundle, say $L$.

Under the splitting $W^\perp \otimes V \cong V \oplus V^*$, the
transposition $\tau$ interchanges $V$ and $V^*$ (and induces the
given identification $V \cong V^*$). It is well-known that the
action of such map $\iota=\tau$ on the Picard group of the
complement $P' \setminus \{0\}$ is non-trivial: we have $\iota^*L
\cong L^{-1}$.

Assume that a small $\iota$-invariant neighborhood $U' \subset P'$
of $0 \subset P'$ admits a smooth projective resolution $X'$
compatible with the symplectic form. By Lemma~\ref{vf.lifts} the map
$f:X' \to U'$ is semismall. In particular, $\codim f^{-1}(0) \subset
X'$ is at least $\frac{1}{2}\dim U' = 2k-1 \geq 2$. Therefore the
Picard group of the variety $X'$ is the same as the Picard group of
the complement $U' \setminus \{0\}$. Thus the line bundle $L$
extends to $X'$, and either $L$ or $L^{-1}$ is an ample line bundle
on $X'$. But $\iota$ interchanges these two bundles. If it were to
extend to an involution $\iota:X' \to X'$, we would get that both
$L$ and $L^{-1}$ are ample line bundles on $X'$. Therefore
$\calo_{X'}$ must be ample, which means that $X'$ is affine, and the
projective map $f:X' \to U'$ from $X'$ to the normal variety $U'$
must be one-to-one. This contradicts the smoothness of $X'$.
\endproof

\begin{remark}
The space $P'$ of course does admit a smooth symplectic resolution,
namely, the total space $T^*P(V)$. The point is that the involution
$\iota$ induces a flop of this resolution, not an automorphism. The
idea of using this involution to prove that the moduli space
$\M(2,0,k)$, $k \geq 3$ does not admit a good resolution is due to
K. O'Grady. We would like to thank him for explaining it to us.
\end{remark}

\proof[Proof of Theorem~\ref{main}.] Assume given a smooth
projective resolution $X \to U_m$ of a small neighborhood $U_m$ of
$0 \subset P(PGL(2),V) \times V$, and assume that the resolution $X$
is compatible with the symplectic form.  Shrinking $U_m$ if
necessary, we can assume that $U_m = U_a \times U_b$, where $U_a
\subset V$, $U_b \subset P = P(PGL(2),V)$ are small neighborhoods of
$0$.

By Lemma~\ref{vf.lifts}, the natural action of the Lie algebra
$\ssp(V)$ on the space $P$ lifts to a $\ssp(V)$-action on the
symplectic resolution $X \to U_b$. Shrinking $U_b$ further, if still
necessary, we can assume that it is $Sp(V)$-invariant. Since $Sp(V)$
is a simply-connected semisimple algebraic group, the action of the
Lie algebra $\ssp(V)$ integrates to an $Sp(V)$-action on the
resolution $X$. In particular, the element $\iota \in Sp(V)$ must
act on the variety $X$.

Consider the open subset $U \subset P$ defined in
Lemma~\ref{dec}. Since $U \subset P$ is dense, the intersection $U
\cap U_b$ is dense in $U_b$. Replacing $U_b$ with $U_b \cap U$,
shrinking it further if necessity persists, and applying
Lemma~\ref{dec}, we can assume that $\rho^{-1}(U_b) = U_c \times
U'$, where $U_c \subset U_1 \subset \wt{P}_1$ is an open subset, and
$U' \subset P'$ is an open $\iota$-invariant neighborhood of $0 \in
P'$. Applying Lemma~\ref{factors}, we conclude that the neighborhood
$U' \subset P'$ must have a $\iota$-equivariant smooth projective
resolution $X'/U'$ compatible with the given symplectic form. This
contradicts Lemma~\ref{no.res}. \endproof

\section{Formality}\label{frm.sec}

Our remaining task is to prove Proposition~\ref{frm}. We will need
one general result on formality in families proved in
\cite{K4}. Namely, recall that for any associative algebra $B$, or,
more generally, for a flat sheaf $\B$ of associative algebras on a
scheme $X$, the {\em Hochschild cohomology sheaves} $\hh^\hdot(\B)$
are defined. Explicitly, they can be computed by the standard complex
with terms
$$
\hhom(\B^{\otimes \hdot},\B)
$$
and a certain differential, whose precise form we will not need. If
the algebra $\B^\hdot$ is graded, the Hochschild cohomology sheaves
inherit the grading; we denote the component of degree $l$ by
$\hh^\hdot_l(\B^\hdot)$.

\begin{lemma}[{{\cite[Theorem 4.2]{K4}}}]\label{Q}
Let $\A^\hdot$ be a DG algebra of flat quasicoherent sheaves on a
reduced irreducible scheme $X$. Let $\B^\hdot$ be the homology
algebra of the DG algebra $\A^\hdot$. Assume that the sheaves
$\B^\hdot$ are coherent and flat on $X$, and that for any integers
$l$, $i$, the degree-$l$ component $\hh^i_l(\B^\hdot)$ of the $i$-th
Hochschild cohomology sheaf $\hh^i(\B^\hdot)$ is also coherent and
flat.
\begin{enumerate}
\item Assume that $X$ is affine. If the fiber $\A^\hdot_x$ is formal
for a generic point $x \in X$, then it is formal for an arbitrary
point $x \in X$.
\item Assume that $\hh^2_l(\B^\hdot)$ has no global sections for all
$l \leq -1$. Then the DG algebra $\A^\hdot_x$ is formal for every
point $x \in X$.\endproof
\end{enumerate}
\end{lemma}

We will apply this result to the study of the ideal sheaf $\I_\z$ of
a $0$-dimensional subscheme $\z \subset S$ of some length $p$ in a
$K3$ surface $S$. We want to prove that the DG algebra
$\RHom^\hdot(\I_\z,\I_\z)$ is formal.

Our first observation is that by Lemma~\ref{Q} we can assume that
$\z \subset S$ is a union of several points with reduced scheme
structure. Indeed, every subscheme $\z \subset S$ can be deformed to
such a subscheme, and the algebra structure on
$\Ext^\hdot(\I_\z,\I_\z)$ does not depend on the choice of $\z$, so
that the Hochschild cohomology sheaves in the assumptions of
Lemma~\ref{Q} are locally trivial. Thus Lemma~\ref{Q}~\thetag{i}
applies and shows that formality is stable under specialization.

Assume therefore that $\z = \bigcup_{1 \leq i \leq p} \{z_i\}
\subset S$ is the union of points. Choose a hyperk\"ahler metric on
the $K3$ surface $S$. Recall that to a hyperk\"ahler metric on a
holomorphically symplectic manifold $M$, one associates a set of
integrable complex structures on $M$, one for each imaginary
quaternion $h$ with $h^2=-1$. The set of such quaternions is
naturally parametrized by the complex projective line $\C P^1$. One
can collect all these complex structures into a single integrable
complex structure on the product $M \times \C P^1$. In this way, one
obtains the so-called holomorphic {\em twistor space } $X$ of the
hyperk\"ahler manifold $M$. The product decomposition $X \cong M
\times \C P^1$ is not holomorphic. However, the natural projection
$\pi:X \to \C P^1$ is holomorphic. Moreover, for every point $m \in
M$ the product $\{ m \} \times \C P^1 \subset X$ is a holomorphic
submanifold, called {\em horizontal section} of the twistor
projection $\pi:X \to \C P^1$ corresponding to $m$. The fiber $X_0$
over the point $0 \in \C P^1$ does not depend on the choice of
hyperk\"ahler metric on $M$, it is canonically isomorphic to the
original complex manifold $M$.

Thus, having chosen a hyperk\"ahler metric on our $K3$ surface $S$,
we obtain the associated twistor deformation $X/\C P^1$ -- the fiber
$X_0$ over $0 \in \C P^1$ is by definition the $K3$ surface $S$
itself, and the fibers over other points in $\C P^1$ correspond to
$S$ with complex structures coming from different quaternions. Lift
each of the points $z_i$ to the corresponding horizontal section
$\ZZ_i \subset X$, $\ZZ_i \cong \C P^1$ of the twistor projection
$\pi:X \to \C P^1$. Let $\I_\ZZ$ be the ideal sheaf of the union
$\bigcup \ZZ_i \subset X$ of these sections on the complex manifold
$X$.

The DG algebra $\RHom^\hdot(\I_\z,\I_\z)$ naturally deforms to a flat
DG algebra $\Rhom^\hdot(\I_\ZZ,\I_\ZZ)$ of coherent sheaves on $\C P^1$
-- to obtain $\Rhom^\hdot(\I_\ZZ,\I_\ZZ)$, one considers the relative
$\Rhom^\hdot$ on $X$ and applies the direct image $R^\hdot\pi_*$
with respect to the twistor projection $\pi:X \to \C P^1$. The
cohomology algebra $\ext^\hdot(\I_\ZZ,\I_\ZZ)$ is a flat deformation
of the algebra $\Ext^\hdot(\I_\z,\I_\z)$.

\begin{defn}
A coherent sheaf on $\C P^1$ is {\em of weight $l$} if it a sum of
several copies of the sheaf $\calo(l)$.
\end{defn}

\begin{lemma}\label{weights}
For every $k \geq 0$, the sheaf $\ext^k(\I_\ZZ,\I_\ZZ)$ on $\C P^1$
is of weight $k$.
\end{lemma}

\proof{} The only $k$ with non-trivial $\ext^k(\I_\ZZ,\I_\ZZ)$ are $k =
0,1,2$. For $k=0$, the sheaf $\ext^0(\I_\ZZ,\I_\ZZ)$ is obviously the
trivial sheaf $\calo = \calo(0)$. For $k=2$, one applies the
relative Serre duality to $X/\C P^1$ and concludes that
$$
\ext^2(\I_\ZZ,K_{X/\C P^1} \otimes \I_\ZZ)^* \cong
\ext^0(\I_\ZZ,\I_\ZZ). 
$$
This proves the claim, since the relative canonical bundle $K_{X/\C
P^1}$ is canonically isomorphic to $\pi^*\calo(2)$ -- the
isomorphism is given by the relative symplectic form $\Omega \in
H^0(X,K_{X/\C P^1})$.

It remains to consider the case $k=1$. It is not diffcult to compute
the sheaf $\ext^1(\I_\ZZ,\I_\ZZ)$ directly. However, we prefer to
use the following geometric argument. The ideal sheaf $\I_\z$
corresponds to a point $z \in S^{[p]}$ in the Hilbert scheme
parametrizing $0$-dimensional subschemes in $S$ of some length $p$
(if $S$ is not algebraic, replace the Hilbert scheme with the Douady
space). Consider the moduli space $X^{[p]}$ of pairs
$$
\left\langle\; \text{a point }I \in \C P^1, \text{ a $0$-dimensional
subscheme }\z \subset X_I \text{ of length } p \;\right\rangle.
$$
The space $X^{[p]}$ projects onto $\C P^1$, and the fiber over a
point $I \in \C P^1$ is the Hilbert scheme $X_I^{[p]}$ (in
particular, the fiber $\left(X^{[p]}\right)_0$ is the Hilbert scheme
$S^{[p]}$). The union $\bigcup \ZZ_i \subset X$ gives a section $Z
\subset \X^{[p]}$, $Z \cong \C P^1$ of the projection $X^{[p]} \to
\C P^1$. By the usual deformation theory, the sheaf
$\ext^1(\I_\ZZ,\I_\ZZ)$ is isomorphic to the normal bundle
$\N_Z(X^{[p]})$.

Now, since the points $z_i \in S$ are distinct by assumption, the
Hilbert scheme $S^{[p]}$ near the point $z \in S^{[p]}$ is
isomorphic to the $p$-fold product $S^p$. Analogously, in a small
neighborhood of $Z \subset X^{[p]}$, the space $X^{[p]}$ is
isomorphic to the $p$-fold product $X^p$ taken over $\C P^1$. But
this product is the twistor space of the product $S^p$, with the
product hyperk\"ahler metric. Therefore, by the general twistor
theory, the normal bundle $\N_Z(X^{[p]})$ is a sum of several copies
of the sheaf $\calo(1)$.  \endproof

Lemma~\ref{weights} allows to finish the proof of
Proposition~\ref{frm} very quickly. 

\proof[Proof of Proposition~\ref{frm}.] Consider the sheaf of DG
algebra $\Rhom^\hdot(\I_\ZZ,\I_\ZZ)$ on the scheme $X = \C P^1$. Its
cohomology sheaf $\B^\hdot = \ext^\hdot(\I_\ZZ,\I_\ZZ)$ is locally
constant as a sheaf of algebras, so that the Hochschild cohomology
sheaves $\hh^\hdot(\B^\hdot)$ are locally trivial, and we can apply
Lemma~\ref{Q}. By Lemma~\ref{weights}, the degree-$l$ component of
the sheaf
$$
\hhom^\hdot(\B^{\hdot \otimes k},\B^\hdot)
$$
is of weight $l$ for any integer $l$ and any $k \geq 0$. But for any
map of sheaves of weight $l$, its kernel and cokernel are obviously
also of weight $l$. Therefore the Hochschild cohomology sheaves
$\hh^k_l(\B^\hdot)$ are also of weight $l$ for any integers $l$ and
$k$. Since sheaves of negative weight have no global sections, the
condition of Lemma~\ref{Q}~\thetag{ii} is also satisfied, and we
conclude that indeed, the DG algebra $\RHom^\hdot(\I_\z,\I_\z)$ is
formal.
\endproof

We conclude this Section with a sort of an extended remark. The idea
behind our proof of formality was discovered by M. Verbitsky
(\cite{ver1}, see also an exposition in \cite{KV}). He considered
deformations of vector bundles on a hyperk\"ahler manifold, and
proved formality by a version of our Lemma~\ref{weights}. We
essentially extend his approach to degenerating sheaves, in the
simplest case of a hyperk\"ahler surface. However, we do need to
assume that the degeneration locus $\z \subset S$ is a reduced
subscheme, and use a roundabout argument to handle the general
case. The reason for this is very simple: for a non-reduced
subscheme $\z \subset S$, the statement of Lemma~\ref{weights} is
{\em false}.

Our proof would work just as well if the deformation $X^{[p]}$ were
the twistor space for some hyperk\"ahler structure on the Hilbert
scheme $S^{[p]}$. This is not the case though. The simplest way to
see this is to note that the fiber $X^{[p]}$ over an arbitrary point
$I \in \C P^1$, being a Hilbert scheme in its own right, contains a
non-trivial divisor -- namely, the locus of non-reduced subschemes
$\z \in X_I$. On the other hand, the generic fiber of the twistor
deformation of an arbitrary hyperk\"ahler manifold has no rational
cohomology classes of Hodge type $(1,1)$.

Algebraically, the sheaf $\ext^1(\I_\ZZ,\I_\ZZ)$ for a non-reduced
scheme $\ZZ \in X^{[p]}$ acquires several summands of type $\calo(2)$
(and the same number of summmands of type $\calo(0)$).

A natural thing to do would be to replace $X^{[p]}$ with the twistor
space for some hyperk\"ahler structure on $S^{[p]}$. However, we do
not have a modular interpretation of these twistor spaces extending
the modular interpretation of the Hilbert scheme $S^{[p]}$.

An analogous problem for the Hilbert scheme of $\C^2$ has been
studied recently in \cite{KKNO}. In that paper the authors managed
to find a modular interpretation of an actual twistor deformation of
the Hilbert scheme of $\C^2$. To do this, however, they had to
consider sheaves on a {\em non-commutative} deformation of the
underlying manifold $\C^2$.

An analogy with \cite{KKNO} suggests that in order to obtain a
modular interpretation of the twistor space for the Hilbert scheme
$S^{[p]}$ -- in particular, in order to obtain a natural proof of
formality -- one has to consider quantizations of the underlying
$K3$ surface $S$.

Of course, the case of a $K3$ surface is much more difficult than
the case of $\C^2$ -- where an explicit quantization is given by the
algebra of differential operators on an affine line. Still, even the
$K3$ case might not be completely beyond modern techniques.
D. Huybrechts informs us that the line $l$ in the period space of
the hyperk\"ahler manifold $S^{[p]}$ defined by the family
$X^{[p]}/\C P^1$, while not being a twistor line, is very close to a
twistor line -- there is a twistor line in an arbitrarily small
neighborhood of $l$. Thus the amount of non-commutativity one has to
introduce into the problem can arbitrarily small. We hope to return
to this in the future.

\subsection*{Acknowledgements.} We would like to thank K. O'Grady
and D. Huybrechts for valuable discussions. Most of the work in this
paper was done while the authors took part in a school on moduli spaces
in Trieste, at ICTP. We would like to thank the organizers for this
opportunity to enjoy the hospitality of the Centre and its wonderful
and stimulating atmosphere.

\noindent
{\sc
Steklov Math Institute\\
Moscow, USSR\\[2mm]
and\\[2mm]
Universit\"at Gutenberg\\
Mainz, Deutschland}

\bigskip

\noindent
{\em E-mail addresses\/}: {\tt kaledin@mccme.ru}\\
\phantom{{\em E-mail addresses\/}: }{\tt lehn@mathematik.uni-mainz.de}

\end{document}